\documentclass[12pt]{article}
\usepackage[english]{babel}
\usepackage{amsmath}
\usepackage{amssymb}
\usepackage{amsthm}
\newcommand{\acks}{\bigskip\par\noindent{\bf Acknowledgement}\par}

    \newtheorem{lemma}{Lemma}[section]
    \newtheorem{proposition}[lemma]{Proposition}
    \newtheorem{theorem}[lemma]{Theorem}
    
     \newtheorem{corollary}[lemma]{Corollary}
    \newtheorem{remark}[lemma]{Remark}

\newcommand{\no}{{\nonumber}}
\newcommand{\ov}{\overline}
\newcommand{\bear}{\begin{array}}
\newcommand{\enar}{\end{array}}
\newcommand{\beq}{\begin{equation}}
\newcommand{\eeq}{\end{equation}}
\newcommand{\beqn}{\begin{eqnarray}}
\newcommand{\eeqn}{\end{eqnarray}}
\newcommand{\beit}{\begin{itemize}}
\newcommand{\eeit}{\end{itemize}}

\renewcommand{\d}{\,{\rm d}}

\newcommand{\rsp}{{\bf R}}
\newcommand{\csp}{{\bf C}}
\newcommand{\nsp}{{\bf N}}

\newcommand{\ds}{\displaystyle}

\newcommand{\ve}{\varepsilon}
\newcommand{\g}{\gamma}
\renewcommand{\l}{\lambda}
\newcommand{\dl}{\delta}

\newcommand{\s}{\sigma}
\renewcommand{\S}{\Sigma}

\newcommand{\Om}{\Omega}
\renewcommand{\a}{\alpha}
\renewcommand{\b}{\beta}

\newcommand{\wtil}[1]{\widetilde{#1}}
\newcommand{\pn}{\par \noindent}
\newcommand{\med}{\medskip}
\newcommand{\qq}{\qquad}
\newcommand{\q}{\quad}
\newcommand{\hookto}{\hookrightarrow}

\title{Optimal time and space regularity for
\\
solutions of degenerate differential equations}
\author{{\rm Alberto Favaron}
\\
\small{Dipartimento di Matematica ``F. Brioschi'' Politecnico di Milano,}
\\
\small{via Bonardi 9, 20133 Milano, Italy. Email: alberto.favaron@polimi.it}
\\
\small{Phone: +39-02-23994639}
\thanks
{The author is a member of the research group GNAMPA of the Italian Istituto
Nazionale di Alta Matematica (INdAM)}}
\date{}

\begin{document}
\maketitle 
\pn
{\bf Abstract.} We derive optimal regularity, in both time and space,
for solutions of the Cauchy problem related to a degenerate differential
equation in a Banach space $X$. Our results exhibit a sort of prevalence
for space regularity, in the sense that the higher is the
order of regularity with respect to space, 
the lower is the corresponding order of regularity with respect to time.
\med\pn
\med\pn

\section{Introduction}\label{Sec1} 
\setcounter{equation}{0}
Let $X$ be a Banach space and let $M$ and $L$ be two closed linear operators from $X$ to itself, whose domains fulfill the relation 
${\cal D}(L)\subset {\cal D}(M)$. 
Further, let $f$ be a continuous function 
from $[0,T]$ into $X$, $T>0$, and let $u_0$ be a given element of $X$.
The question of maximal regularity for the initial value problem
\beqn\label{1.1}
&&\hskip -1truecm
\left\{\!\!
\begin{array}{lll}
D_t(Mv(t))=Lv(t)+f(t),\q t\in(0,T],\qq\Big(D_t=\ds\frac{d}{dt}\Big)
\\[2mm]
Mv(0)=u_0,
\end{array}
\right.
\eeqn
concerns what kind of properties, in time and/or in space, 
the data need satisfy, in order that the solution $v$ to (\ref{1.1}) 
exists and the derivative $D_tMv$ possesses similar regularity as the data.

Since the natural operator associated to (\ref{1.1}) is $A=LM^{-1}$, 
we are led to consider the equivalent problem
\beqn\label{1.2}
&&\hskip -2truecm
\left\{\!\!
\begin{array}{lll}
D_tw(t)=Aw(t)+f(t),\q t\in(0,T],
\\[2mm]
w(0)=u_0,
\end{array}
\right.
\eeqn
where $w=Mv$. Hence, the question of maximal regularity for (\ref{1.1}) 
is strictly related to the regularity of the semigroup generated by $A$.
This yields to the analysis of the spectral equation $\l u-Au=f$, $\l\in\csp$, $f\in X$, 
in order to obtain an estimate of type
\beqn\label{1.3}
&&\hskip -1,5truecm
\|(\l I-A)^{-1}f\|_X\le C(|\l|+\l_0)^{-\b}\|f\|_X, \qq\l\in\S,
\eeqn
$I$ being the identity operator. Here, $\b\in(0,1]$, $\l_0\ge 0$ is large enough, 
and $\S$ is a complex region containing the half plane $\Re{\rm e}\l\ge 0$.

Of course (cf. \cite[Theorem 3.17]{FY1} and \cite[Theorem 1]{FY2}),  
if $A$ satisfies assumption (\ref{1.3}) with $\b=1$, then the results of 
maximal regularity are analogous to those exhibited in \cite{SI}
for the non degenerate case, corresponding to $M=I$ in (\ref{1.1}), 
and for which, nowadays, a wide literature exists. 
In particular, in the case $\b=1$, 
$D_tMv$ has exactly the same regularity as the data.
This extends to (\ref{1.1}) the known results on the maximal regularity of
solutions to (\ref{1.2}) when the semigroup generated by $A$ is analytic.
On the contrary, according to \cite{FY1} and \cite{HP},
if $\b\in(0,1)$, then, in general, the semigroup
generated by $A$ is no longer analytic, but only infinitely differentiable,
and this implies that the time regularity of the solutions to (\ref{1.1})
decreases. We refer to \cite[Theorem 3.26]{FY1}, \cite[Theorem 9]{FY2} 
and \cite[Theorem 7.2]{FLT} for precise statements 
and amounts of the loss of regularity, but, briefly,
the quoted theorems say that if $f\in C^\tau([0,T];X)$, $\tau\in(1-\b,1)$, 
and $u_0$ fulfills some natural consistency conditions, 
then $D_tMv\in C^\nu([0,T];X)$, where $\nu=\tau+\b-1$.

Notice that, at present, one of the main deficiencies in the theory of degenerate 
equations is the absence of results of space regularity, in the case
when $\b\in(0,1)$ in (\ref{1.3}). It is our aim, here, to give a contribution
in this field providing an optimal ``cross" regularity result, in which
both time and space regularity for $D_tMv$ are established. As we shall see,
the space regularity prevails, in the sense that the increase in space regularity
reflects in a decrease of the order of time regularity.

The plan of the paper is the following.
In Section \ref{Sec2}, to a linear operator $A$ from a
Banach space $X$ to itself, whose resolvent satisfies (\ref{1.3}) 
in a region $\S$ depending on a additional parameter $\a\in(\b,1]$, 
we associate the corresponding infinitely differentiable semigroup 
$\{{\rm e}^{tA}\}_{t\ge 0}$ on $X$. Moreover, we recall the definition
of the interpolation spaces $(X,{\cal D}(A))_{\g,p}$ between 
the domain ${\cal D}(A)$ of $A$ and $X$. 

Section \ref{Sec3} is devoted to show that the uniform norm 
$\|A^n{\rm e}^{tA}\|_{{\cal L}(X;(X,{\cal D}(A))_{\g,p})}$ blows up, 
as $t$ goes to $0$, as a suitable negative power of $t$ depending 
on $\a$, $\b$, $\g$ and $n$. 
Further, the blow-up rate is greater than the one observed in \cite{LU} 
for the non degenerate case. As a corollary, we show that 
for every $\ve\in(0,T]$ and $\s\in(0,1)$ the map $t\to A^n{\rm e}^{tA}$ belongs to $C^\s([\ve,T];{\cal L}(X;(X,{\cal D}(A))_{\g,p}))$.

Using the results of Section \ref{Sec3}, in Section \ref{Sec4} we
establish time and space regularity of some basic operator functions,
which appear naturally when $Mv$ and $D_tMv$
are represented in terms of the Volterra integral equation 
equivalent to (\ref{1.2}). In particular, Lemmas \ref{lem4.1}--\ref{lem4.5}
highlight the above mentioned fact that the higher is the order $\g$ 
of the interpolation space $(X,{\cal D}(A))_{\g,p}$ where we look
for space regularity, the lower is the H\"older exponent $\s$ of regularity in time.

Section \ref{Sec5} contains our main results.
First, using Lemmas \ref{lem4.1}--\ref{lem4.3}, 
in Theorems \ref{thm5.1} and \ref{thm5.3} 
we show that if $\g$ and $\s$ are opportunely chosen, $\s<\g$, 
then (\ref{1.1}) has a unique strict solution $v$ such that 
$Mv\in C^{\s}([0,T];(X,{\cal D}(A))_{\g,p})$. 
Then, combining Lemmas \ref{lem4.3}--\ref{lem4.5}, in Theorem \ref{thm5.4}
we prove that, if $\a$ and $\b$ are large enough and the data pair 
$(f,u_0)$ satisfies some suitable space--time assumptions, 
the regularity $C^{\s}([0,T];(X,{\cal D}(A))_{\g,p})$ holds
for the derivative $D_tMv$, too.
\section{Preliminary material and notations}\label{Sec2}
\setcounter{equation}{0}
Let $X$ be a Banach space endowed with norm $\|\cdot\|_X$ 
and let $A:{\cal D}(A)\subset X\to X$ be a single valued linear operator.
Recalling that the resolvent set $\rho(A)$ of $A$ is the set of values
$\l\in \csp$ such that $\l I-A$ 
has a bounded inverse $(\l I-A)^{-1}$ with domain dense in $X$,
we assume that $A$ satisfies the
following resolvent condition: 
\begin{itemize}
\item[(H1)] $\rho(A)$ contains the complex region
$\S=\{\l\in\csp:\Re{\rm e}\l\ge -c(|\Im{\rm m}\l|+1)^\a\}$ and,
for every $\l\in\S$, the following estimate holds
\beqn
&&\hskip -2truecm
\|(\l I-A)^{-1}\|_{{\cal L}(X)}\le C(|\l|+1)^{-\b},\no
\eeqn 
for some exponents $0<\b<\a\le 1$ and constants $c$, $C>0$.
\end{itemize}
Here, as usual, ${\cal L}(X)$ denotes  
the Banach space ${\cal L}(X;X)$ of all bounded linear operators from $X$ to $X$, 
equipped with the uniform operator norm.

According to \cite[Theorem 3.1]{FY1}, assumption (H1) implies
that $A$ generates an infinitely differentiable semigroup on $X$. 
More precisely, introduce the family 
$\{{\rm e}^{tA}\}_{t>0}\subset{\cal L}(X)$ defined by the Dunford integral
\beqn\label{2.1}
&&\hskip -2truecm
{\rm e}^{tA}={(2\pi i)^{-1}}\int_{\Gamma}{\rm e}^{t\l}(\l I-A)^{-1}\d\l,\q\ t>0,
\eeqn
where $\Gamma\subset\S$ is the contour parametrized by 
$\l=-c(|\eta|+1)^\a+i\eta$, $-\infty<\eta<\infty$. Define also ${\rm e}^{0A}=1$.
Then $\{{\rm e}^{tA}\}_{t\ge 0}$ is a semigroup on $X$, 
infinitely many times differentiable for $t>0$ with 
$D_t{\rm e}^{tA}=A{\rm e}^{tA}$. In addition,
${\rm e}^{tA}$ satisfy the estimates (see \cite[Proposition 3.2]{FY1})
\beqn\label{2.2}
&&\hskip -1truecm
\|A^k{\rm e}^{tA}\|_{{\cal L}(X)}\le \wtil c_kt^{(\b-k-1)/\a},\q t>0,\q k\in\nsp\cup\{0\},
\eeqn
where the $\wtil c_k$'s are positive constants depending on $k$. 
Of course, due to (\ref{2.2}), if $\b<1$, then the function $t\to{\rm e}^{tA}$ 
is not bounded as $t\to 0^+$.
As a consequence, ${\rm e}^{tA}$ is not necessarily strongly continuous
in the norm of $X$ on the subspace $\ov{{\cal D}(A)}$.

We stress that, 
even though here we are following the approach in \cite{FY1}, 
resolvent conditions of type (H1) were already introduced in \cite{HP}. 
In particular, in \cite[Remark page 383]{HP} it was showed 
that if $U$ is a closed linear operator with dense domain
and such that
\beqn\no
&&\hskip 0,5truecm
\|(\l I-U)^{-1}\|_{{\cal L}(X)}\le C(\Re{\rm e}\l +|\Im{\rm m}\l|^\b)^{-1},\q\;
\Re{\rm e}\l>0,\q\b\in(0,1),
\eeqn
then $U$ generates a semigroup ${\rm e}^{tU}$ which is infinitely differentiable for $t>0$.

For our purposes, we recall now the definitions 
of two classes of real interpolation spaces between
${\cal D}(A)$ and $X$.
First of all, we specify a topology on ${\cal D}(A)$ equipping it with the norm
$\|y\|_{{\cal D}(A)}=\|y\|_X+\|Ay\|_X$ which makes ${\cal D}(A)$ a Banach space.
Now, if $Z$ is a Banach space, for an $Z$-valued strongly measurable function
$g(\xi)$, $\xi\in(0,\infty)$, we set
\beqn\no
&&\hskip -1truecm
\|g\|_{L_*^p(Z)}=\Big(\int_0^\infty\|g(\xi)\|_Z^p\frac{\d\xi}{\xi}\Big)^{1/p},\qq p\in[1,\infty),
\\[2mm]
&&\hskip -1truecm
\|g\|_{L_*^\infty(Z)}=\sup_{\xi\in(0,\infty)}\|g(\xi)\|_Z,\qq p=\infty.\no
\eeqn
Then, according to \cite[pag. 26]{FY1}, 
for every $\g\in(0,1)$ and $p\in[1,\infty]$ we introduce
the (intermediate) spaces 
\beqn\no
&&\hskip -1truecm
X_A^{\g,p}=\{x\in X: \|\xi^\g A(\xi I-A)^{-1}x\|_{L_*^p(X)}<\infty\},
\eeqn
which becomes Banach spaces when endowed with the norm
\beqn\no
&&\hskip -1,5truecm
\|x\|_{X_A^{\g,p}}:=\|x\|_X+ \|\xi^\g A(\xi I-A)^{-1}x\|_{L_*^p(X)}.
\eeqn
Also, for $\g\in(0,1)$ and $p\in[1,\infty]$ we denote with 
$V(p,\g,{\cal D}(A),X)$ the space of all $X$-valued functions
$v(\xi)$, $\xi\in(0,\infty)$, having the property that the maps
$\xi\to\xi^\g v(\xi)$ and $\xi\to\xi^\g v'(\xi)$
belong, respectively, to $L_*^p({\cal D}(A))$ and $L_*^p(X)$.
As it is well-known (cf. \cite[Lemma 1.8.1]{TR}), 
the spaces $V(p,\g,{\cal D}(A),X)$ are Banach 
spaces with the norm
\beqn\no
&&\hskip -2truecm
\|v\|_{V(p,\g,{\cal D}(A),X)}=
\|\xi^\g v\|_{L_*^p({\cal D}(A))}+\|\xi^\g v'\|_{L_*^p(X)}
\eeqn
and any function $v\in V(p,\g,{\cal D}(A),X)$ has a $X$-valued continuous extension 
at $t=0$. This lead to define the trace spaces (cf. \cite[Theorem 1.8.2]{TR})
\beqn\no
&&\hskip -2truecm
(X,{\cal D}(A))_{\g,p}=\{x\in X: x=v(0),\ v\in V(p,1-\g,{\cal D}(A),X)\},
\eeqn
which turn out to be real interpolation spaces between
${\cal D}(A)$ and $X$. They are Banach spaces endowed with the norm
\beqn\no
&&\hskip -1truecm
\|x\|_{(X,{\cal D}(A))_{\g,p}}=
\inf\{\|v\|_{V(p,1-\g,{\cal D}(A),X)}: x=v(0), \ v\in V(p,1-\g,{\cal D}(A),X)\}.
\eeqn
Further, for $\g\in(0,1)$ and $1\le p_1<p_2\le\infty$ we have
\beqn\no
&&\hskip -1,5truecm
{\cal D}(A)\hookto (X,{\cal D}(A))_{\g,p_1}\hookto
(X,{\cal D}(A))_{\g,p_2}\hookto\ov{{\cal D}(A)},
\eeqn
whereas, for $0<\g_1<\g_2<1$, we have
\beqn\no
&&\hskip -1,5truecm
(X,{\cal D}(A))_{\g_2,\infty}\hookto (X,{\cal D}(A))_{\g_1,1}.
\eeqn
The classes $X_A^{\g,p}$ and $(X,{\cal D}(A))_{\g,p}$ are related one to each other
by the following continuous embedding (see \cite[Theorem 2]{WI})
\beqn\label{2.32.3}
&&\hskip -1truecm
X_A^{\g,p}\hookto(X,{\cal D}(A))_{\g,p},\q \g\in(0,1),\q p\in[1,\infty],
\\[1mm]
&&\hskip -1truecm
(X,{\cal D}(A))_{\g,p}\hookto X_A^{\g+\b-1,p},\q\g\in(1-\b,1),\q p\in[1,\infty],\no
\eeqn
which become identities with equivalence of the respective norms when $\b=1$ in (H1).

As we said before in Introduction, the natural
operator $A$ associated to (\ref{1.1}) is the operator $LM^{-1}$ having
domain ${\cal D}(A)=M({\cal D}(L))$. From this point of view, 
in despite of the non degenerate case where the characterizations are wide, 
in the degenerate case a characterization of either $X_A^{\g,p}$ 
or $(X,{\cal D}(A))_{\g,p}$ is still lacking, even in the very common situation where $X=L^p(\Om)$, $p\in[1,\infty]$, $\Om\subset\rsp^n$, $M$ is the multiplication operator by a fixed positive function $m\in L^{\infty}(\Om)$ which may have zeros,
and $L$ is a elliptic second-order linear differential operator with domain
${\cal D}(L)=W^{2,p}(\Om)\cap W_0^{1,p}(\Om)$.  
Up to now, to the author's knowledge, the only available 
result in this direction is Lemma 3.2 in \cite{LT}
where $(X,{\cal D}(A))_{\g,p}$ is shown to contain
a special class of fractional Sobolev spaces. 
Therefore, since in the concrete situations 
the spaces $(X,{\cal D}(A))_{\g,p}$ seem easier to handle than the $X_A^{\g,p}$, 
here we prefer to invoke this class of intermediate spaces 
when interpolation is needed from Proposition \ref{prop3.1} onward.
This choice, opposite to that of \cite{FY1}, seems however to be fruitful
and allows us to improve some estimates for the semigroup ${\rm e}^{tA}$ 
which are known only in the $X_A^{\g,\infty}$-setting 
(see the following Remark \ref{rem3.2}).

We conclude the section introducing some notations we will largely use in the sequel.
Given a Banach space $Z$, $C([0,T];Z)$ and $C^\dl([0,T];Z)$, $\dl\in(0,1)$, 
denote, respectively, the spaces of all continuous
and $\dl$-H\"older continuous functions from $[0,T]$ into $Z$.
The shortenings $\|\cdot\|_{0,T;Z}$ and $\|\cdot\|_{\dl,T;Z}$ stand, respectively,
for the usual sup-norm of $C([0,T];Z)$ and the norm $\|\cdot\|_{C^\dl([0,T];Z)}$
of $C^\dl([0,T];Z)$, i.e.
\beqn\no
&&\hskip -0,5truecm
\|f\|_{0,T;Z}=\sup_{t\in[0,T]}\|f(t)\|_Z,
\\
&&\hskip -0,5truecm
\|f\|_{\dl,T;Z}=\|f\|_{0,T;Z}+|f|_{\dl,T;Z},\q
|f|_{\dl,T;Z}:=\sup_{0\le s<t\le T}\frac{\|f(t)-f(s)\|_Z}{(t-s)^{\dl}}.\no
\eeqn 
Moreover, $B([0,T];Z)$ and $C^1((0,T];Z)$ denote, respectively,
the space of all bounded functions from $[0,T]$ into $Z$ with the sup-norm,
and the space of all strongly differentiable functions on $(0,T]$ 
whose derivatives are continuous from $(0,T]$ into $Z$.
Finally, if $Z_1$ and $Z_2$ are two different Banach space, ${\cal L}(Z_1,Z_2)$
is the Banach space of bounded linear operators from $Z_1$ into $Z_2$ with the usual
uniform operator norm.
\section{Regularity of ${\rm e}^{tA}$ 
with respect the spaces $(X,{\cal D}(A))_{\g,p}$}\label{Sec3}
\setcounter{equation}{0}
Here we show two preliminary results
concerning the behaviour of $A^n{\rm e}^{tA}$, $n\in\nsp\cup\{0\}$,
with respect the interpolation spaces $(X,{\cal D}(A))_{\g,p}$. Essentially, the following
Proposition \ref{prop3.1} says that, when $t$ goes to zero, the norm
$\|A^n{\rm e}^{tA}\|_{{\cal L}(X;(X,{\cal D}(A))_{\g,p})}$ may goes to infinity, 
but not faster than a precise negative power of $t$ depending on $n$, $\g$ and
the exponents $\a$, $\b$ appearing in (H1). 
A similar result is shown in \cite[Proposition 2.3.9]{LU} 
for the non degenerate case, and in \cite[Proposition 3.2]{FY1} 
for the degenerate one. However, in \cite{FY1}, only the case $n=0$ is treated 
and the role of the spaces $(X,{\cal D}(A))_{\g,p}$ is played there 
by the spaces $X_A^{\g,\infty}$.
\begin{proposition}\label{prop3.1}
Let $\a$, $\b\in(0,1]$, $\b<\a$, $\g\in(0,1)$, $p\in[1,\infty]$ and $n\in\nsp\cup\{0\}$.
Then there exist positive constants $C=C(\g,p,n)$ and $C'=C'(\a,\b,\g,p,n)$
such that the following estimates hold
\beqn\label{3.1}
&&\hskip -2truecm
\left\{\!
\begin{array}{lll}
(i)\q\; \|A^n{\rm e}^{tA}\|_{{\cal L}(X;(X,{\cal D}(A))_{\g,p})}
\le Ct^{(\b-n-1-\g)/\a}\,,\q t\in(0,1],
\\[2mm]
(ii)\q\|A^n{\rm e}^{tA}\|_{{\cal L}(X;(X,{\cal D}(A))_{\g,p})}
\le C't^{(\b-n-1)/\a}\,,\q t\ge 1.
\end{array}
\right.
\eeqn
In particular, setting $c_1(T)=C+C'T^{\g/\a}$, for every $T>0$ we obtain
\beqn\label{3.2}
&&\hskip -2truecm
\|A^n{\rm e}^{tA}\|_{{\cal L}(X;(X,{\cal D}(A))_{\g,p})}
\le c_1(T)t^{(\b-n-1-\g)/\a},\q\forall\,t\in (0,T].
\eeqn
\end{proposition}
\begin{proof}
First, for $t\in(0,1]$ and $x\in X$, using the interpolation inequality
$\|y\|_{(X,{\cal D}(A))_{\g,p}}\le c(\g,p)\|y\|_X^{1-\g}\|y\|_{{\cal D}(A)}^{\g}$,
$y\in{\cal D}(A)$, and the estimate 
$\|A^k{\rm e}^{tA}\|_{{\cal L}(X)}\le \wtil c_kt^{(\b-k-1)/\a}$, $t>0$, $k\in\nsp\cup\{0\}$, 
we get
\beqn
&&\hskip -1truecm
\|t^nA^n{\rm e}^{tA}x\|_{(X,{\cal D}(A))_{\g,p}}
\le c(\g,p)\|t^nA^n{\rm e}^{tA}x\|_X^{1-\g}\|t^nA^n{\rm e}^{tA}x\|_{{\cal D}(A)}^\g\no
\\[2mm]
&&\hskip -1,5truecm
\le c(\g,p)\wtil c_n^{\,1-\g}t^{(1-\g)[\b+(\a-1)n-1]/\a}\|x\|_X^{1-\g}
(\|t^nA^n{\rm e}^{tA}x\|_X+\|t^nA^{n+1}{\rm e}^{tA}x\|_X)^{\g}\no
\\[2mm]
&&\hskip -1,5truecm
\le c(\g,p)\wtil c_n^{\,1-\g}t^{(1-\g)[\b+(\a-1)n-1]/\a}
(\wtil c_nt^{[\b+(\a-1)n-1]/\a}+\wtil c_{n+1}t^{[\b+(\a-1)n-2]/\a})^\g\|x\|_X\no
\\[2mm]
&&\hskip -1,5truecm
\le c(\g,p)\wtil c_n^{\,1-\g}(\wtil c_n+\wtil c_{n+1})^\g t^{(1-\g)[\b+(\a-1)n-1]/\a}
t^{\g[\b+(\a-1)n-2]/\a}\|x\|_X\no
\\[2mm]
&&\hskip -1,5truecm
\le c(\g,p){\wtil c}_n^{\,1-\g}(\wtil c_n+\wtil c_{n+1})^\g t^{[\b+(\a-1)n-1-\g]/\a}\|x\|_X.\no
\eeqn
This proves (\ref{3.1})(i) with 
$C= c(\g,p){\wtil c}_n^{\,1-\g}(\wtil c_n+\wtil c_{n+1})^\g$.
Concerning (\ref{3.1})(ii), instead, for $t\ge 1$ and $x\in X$, 
using (\ref{3.1})(i) with $t=1/2$ and $n=0$, we easily derive:
\beqn
&&\hskip-1,1truecm
\|t^nA^n{\rm e}^{tA}x\|_{(X,{\cal D}(A))_{\g,p}}
=\Big\|\Big(\frac{t}{t-1/2}\Big)^n(t-1/2)^nA^n
{\rm e}^{(t-\frac{1}{2})A}{\rm e}^{\frac{1}{2}A}x\Big\|_{(X,{\cal D}(A))_{\g,p}}\no
\\[2mm]
&&\hskip-1,5truecm
\le\Big(\frac{t}{t-1/2}\Big)^n\wtil c_n(t-1/2)^{[\b+(\a-1)n-1]/\a}
\|{\rm e}^{\frac{1}{2}A}\|_{{\cal L}(X;(X,{\cal D}(A))_{\g,p})}\|x\|_X\no
\\[2mm]
&&\hskip-1,5truecm
\le 2^n\wtil c_n[1-1/(2t)]^{[\b+(\a-1)n-1]/\a}t^{[\b+(\a-1)n-1]/\a}
\|{\rm e}^{\frac{1}{2}A}\|_{{\cal L}(X;(X,{\cal D}(A))_{\g,p})}\|x\|_X\no
\\[2mm]
&&\hskip-1,5truecm
\le 2^{(n+1-\b)/\a}\wtil c_nt^{[\b+(\a-1)n-1]/\a}
\|{\rm e}^{\frac{1}{2}A}\|_{{\cal L}(X;(X,{\cal D}(A))_{\g,p})}\|x\|_X\no
\\[2mm]
&&\hskip-1,5truecm
\le 2^{[2(n+1-\b)+\g]/\a}\wtil c_nCt^{[\b+(\a-1)n-1]/\a}\|x\|_X\,.\no
\eeqn
This proves (\ref{3.1})(ii) with $C'=2^{[2(n+1-\b)+\g]/\a}\wtil c_nC$. 
Hence, if $1\le T$, from (\ref{3.1})(ii) we get
$\|A^n{\rm e}^{tA}\|_{{\cal L}(X;(X,{\cal D}(A))_{\g,p})}
\le C'T^{\g/\a}t^{(\b-n-1-\g)/\a}$, $1\le t\le T$.
Combining this latter inequality with (\ref{3.1})(i) we obtain (\ref{3.2}).
\end{proof}
\begin{remark}\label{rem3.2}
\emph{ Observe that, due to the continuous embedding (\ref{2.32.3}), 
if $(n,p)=(0,\infty)$ then (\ref{3.2}) agrees with estimate 
$\|{\rm e}^{tA}\|_{{\cal L}(X;X_A^{\g,\infty})}\le C_\g t^{(\b-1-\g)/\a}$ 
in \cite[Proposition 3.2]{FY1}. However,
since $(X,{\cal D}(A))_{\g,p}\hookto (X,{\cal D}(A))_{\g,\infty}$ for $p\in[1,\infty)$,
our estimate really refines that in \cite{FY1}, even in the case $n=0$.}
\end{remark}
Proposition \ref{prop3.1} easily implies
that, when $t$ is bounded away from zero, then,
for every $\s\in(0,1)$, the operator function $t\to A^n{\rm e}^{tA}$ 
is $\s$-H\"older continuous in time with values in ${\cal L}(X;(X,{\cal D}(A))_{\g,p})$.
Indeed, the following corollary holds.
\begin{corollary}\label{cor3.3}
Let $\a$, $\b\in(0,1]$, $\b<\a$, $\g\in(0,1)$, $p\in[1,\infty]$ and $n\in\nsp\cup\{0\}$.
Then, for every $\s\in(0,1)$ and $0<s<t\le T$ we have
\beqn\label{3.3}
&&\hskip -1,5truecm
\|A^n{\rm e}^{tA}-A^n{\rm e}^{sA}\|_{{\cal L}(X;(X,{\cal D}(A))_{\g,p})}
\le  \s^{-1}c_1(T)s^{(\a+\b-n-2-\g-\a\s)/\a}(t-s)^{\s}.
\eeqn
\end{corollary}
\begin{proof}
For every $x\in X$ and $0<s<t\le T$,
using the identity $[A^n{\rm e}^{tA}-A^n{\rm e}^{sA}]x=
\int_{s}^{t}A^{n+1}{\rm e}^{rA}x\d r$, 
inequality (\ref{3.2}) with $n$ replaced by $n+1$ and
the well-known inequality
$t^\g-s^\g\le(t-s)^\g$, $\g\in(0,1)$, we easily obtain
\beqn\no
&&\hskip -2truecm
\|[A^n{\rm e}^{tA}-A^n{\rm e}^{sA}]x\|_{(X,{\cal D}(A))_{\g,p}}
\le c_1(T)\|x\|_X\int_{s}^{t}\xi^{(\b-n-2-\g)/\a}\d\xi\no
\\[1mm]
&&\hskip 2,9truecm
\le c_1(T)\|x\|_Xs^{(\a+\b-n-2-\g-\a\s)/\a}\int_{s}^{t}\xi^{\s-1}\d\xi\no
\\[1mm]
&&\hskip 2,9truecm
\le \s^{-1}c_1(T)\|x\|_Xs^{(\a+\b-n-2-\g-\a\s)/\a}(t-s)^\s.\no
\eeqn
This completes the proof.
\end{proof}
\section{Time and space regularity of the basic operator functions}\label{Sec4}
\setcounter{equation}{0}
Proposition \ref{prop3.1} and Corollary \ref{cor3.3} enable us to prove
some H\"older-in-time regularity with respect the spaces $(X,{\cal D}(A))_{\g,p}$
for those operator functions involving the semigroup ${\rm e}^{tA}$
which we will encounter later.  
Through the rest of the paper, $c_j(T)$, $j= 2, 3,\ldots$, shall denote 
positive nondecreasing functions of $T$ depending also on $\a$, $\b$, $\g$, $p$
and $\s\in(0,1)$.
\begin{lemma}\label{lem4.1}
Let $\a$, $\b\in(0,1]$ such that $\b<\a$ and $2\a+\b>2$. Then, for every
$\g\in(0,2\a+\b-2)$ and $\s\in(0,(2\a+\b-2-\g)/\a)$ the linear operator
\beqn\label{4.1}
&&\hskip -2truecm
[Q_1g](t):=\int_0^t{\rm e}^{(t-\xi)A}g(\xi)\d\xi
\eeqn
maps $C([0,T];X)$ into $C^{\s}([0,T];(X,{\cal D}(A))_{\g,p})$, $p\in[1,\infty]$,
and satisfies the estimate:
\beqn\label{4.2}
&&\hskip -1truecm
\|Q_1g\|_{\s,T;(X,{\cal D}(A))_{\g,p}}
\le T^{(2\a+\b-2-\g-\a\s)/\a}c_2(T)\|g\|_{0,T;X}.
\eeqn
\end{lemma}
\begin{proof}
First, for every $t\in[0,T]$, inequality (\ref{3.2}) with $n=0$ implies
\beqn\label{4.3}
&&\hskip -1truecm
\|[Q_1g](t)\|_{(X,{\cal D}(A))_{\g,p}}\le
\frac{\a c_1(T)}{\a+\b-1-\g}\|g\|_{0,T;X}t^{(\a+\b-1-\g)/\a},
\eeqn
where the exponent $(\a+\b-1-\g)/\a$ is positive since $0<\g<2\a+\b-2\le\a+\b-1$.
Moreover, when $\s\in(0,(2\a+\b-2-\g)/\a)$ and $0< s<t\le T$
\footnote{Since $\s\in(0,(2\a+\b-2-\g)/\a)$, the case $s=0$ follows from inequality (\ref{4.3}) once we observe that 
$$(\a+\b-1-\g)/\a-\s>\big(\a+\b-1-\g-(2\a+\b-2-\g)\big)/\a=(1-\a)/\a\ge 0$$},
from both (\ref{3.2}) and (\ref{3.3}) with $n=0$ we obtain
\beqn\label{4.4}
&&\hskip -0,5truecm
\|[Q_1g](t)-[Q_1g](s)\|_{(X,{\cal D}(A))_{\g,p}}\no
\\[1mm]
&&\hskip -1truecm
\le c_1(T)\|g\|_{0,T;X}
\Big[\s^{-1}(t-s)^\s\int_0^{s}(s-\xi)^{(\a+\b-2-\g-\a\s)/\a}\d\xi\no
\\[1mm]
&&\hskip 2,2truecm
+\int_{s}^{t} (t-\xi)^{(\b-1-\g)/\a}\d\xi\Big]\no
\\[2mm]
&&\hskip -1truecm
\le \a c_1(T)\|g\|_{0,T;X}
\Big[\frac{\s^{-1} s^{(2\a+\b-2-\g-\a\s)/\a}(t-s)^\s}{2\a+\b-2-\g-\a\s}
+ \frac{(t-s)^{(\a+\b-1-\g)/\a}}{\a+\b-1-\g}\Big]\no
\\[2mm]
&&\hskip -1truecm
\le \a c_1(T)\|g\|_{0,T;X}
\Big[\frac{\s^{-1} s^{(2\a+\b-2-\g-\a\s)/\a}}{2\a+\b-2-\g-\a\s}
+\frac{(t-s)^{(\a+\b-1-\g-\a\s)/\a}}{\a+\b-1-\g}\Big](t-s)^\s.
\eeqn
Finally, summing up (\ref{4.3}) and (\ref{4.4}), we derive (\ref{4.2})
with 
\beqn\label{4.5}
&&\hskip -0,5truecm
c_2(T)=\a c_1(T)\Big[\frac{T^{(1-\a)/\a}(T^\s+1)}{\a+\b-1-\g}
+\frac{\s^{-1}}{2\a+\b-2-\g-\a\s}\Big].
\eeqn
This completes the proof.
\end{proof}
If $g$ is not only merely continuos from $[0,T]$ to $X$, but 
$\s$-H\"older continuous, then the thesis 
of Lemma \ref{lem4.1} follows by a weaker assumption on $\a$ and $\b$
and for larger values of $\g$ and $\s$. Indeed, the proof can be modified in order to avoid Corollary \ref{cor3.3} as it is shown in the following lemma.
\begin{lemma}\label{lem4.2}
Let $\a$, $\b\in(0,1]$ such that $\b<\a$ and $\a+\b>1$. Then, for every
$\g\in(0,\a+\b-1)$ and $\s\in(0,(\a+\b-1-\g)/\a)$ the linear operator $Q_1$
defined by $(\ref{4.1})$ maps $C^\s([0,T];X)$ into $C^{\s}([0,T];(X,{\cal D}(A))_{\g,p})$, $p\in[1,\infty]$, and satisfies the estimate:
\beqn\label{4.6}
&&\hskip -1truecm
\|Q_1g\|_{\s,T;(X,{\cal D}(A))_{\g,p}}
\le T^{(\a+\b-1-\g-\a\s)/\a}c_3(T)\|g\|_{\s,T;X}.
\eeqn
If $g$ is a constant function then 
$(\ref{4.6})$ can be improved until the value $\s=(\a+\b-1-\g)/\a$.
\end{lemma}
\begin{proof}
Since $g\in C^{\s}([0,T];X)$, $\s\in(0,(\a+\b-1-\g)/\a)$, 
when $0< s<t\le T$
\footnote{Since $\s\in(0,(\a+\b-1-\g)/\a)$, the case $s=0$ follows from inequality (\ref{4.3}).} from (\ref{3.2}) with $n=0$ it follows
\beqn\label{4.7}
&&\hskip -0,5truecm
\|[Q_1g](t)-[Q_1g](s)\|_{(X,{\cal D}(A))_{\g,p}}\no
\\[1mm]
&&\hskip -1truecm
\le\int_0^{s}\|{\rm e}^{\xi A}[g(t-\xi)-g(s-\xi)]\|_{(X,{\cal D}(A))_{\g,p}}\d\xi
+\int_{s}^{t} \|{\rm e}^{\xi A}g(t-\xi)\|_{(X,{\cal D}(A))_{\g,p}}\d\xi\no
\\[1mm]
&&\hskip -1truecm
\le c_1(T)\|g\|_{\s,T;X}
\Big[(t-s)^{\s}\int_0^{s}\xi^{(\b-1-\g)/\a}\d\xi+\int_{s}^{t}\xi^{(\b-1-\g)/\a}\d\xi\Big]\no
\\[1mm]
&&\hskip -1truecm
\le \frac{\a c_1(T)}{\a+\b-1-\g}\|g\|_{\s,T;X}[(t-s)^\s s^{(\a+\b-1-\g)/\a}
+(t-s)^{(\a+\b-1-\g)/\a}]\no
\\[2mm]
&&\hskip -1truecm
\le \frac{\a c_1(T)}{\a+\b-1-\g}\|g\|_{\s,T;X}[s^{(\a+\b-1-\g)/\a}
+(t-s)^{(\a+\b-1-\g-\a\s)/\a}](t-s)^\s.
\eeqn
Summing up (\ref{4.3}) and (\ref{4.7}), we derive (\ref{4.6})
with 
\beqn\label{4.8}
&&\hskip -1truecm
c_3(T)=\a(\a+\b-1-\g)^{-1} c_1(T)[2T^{\s}+1].
\eeqn
Last assertion trivially follows simply observing 
that $g(t-\xi)-g(s-\xi)=0$ in the previous computations.
\end{proof}
\begin{lemma}\label{lem4.3}
Let $\a$, $\b\in(0,1]$ such that $\b<\a$ and $\a+\b>1$. 
Moreover, let $x\in{\cal D}(A)$, $\g\in (0,\a+\b-1)$ and  $\s\in(0,(\a+\b-1-\g)/\a)$. 
Then ${\rm e}^{\cdot A}x\in C^{\s}([0,T];(X,{\cal D}(A))_{\g,p})$, $p\in[1,\infty]$,
and satisfies the estimate:
\beqn\label{4.9}
&&\hskip -1,5truecm
\|{\rm e}^{\cdot A}x\|_{\s,T;(X,{\cal D}(A))_{\g,p}}\le c_4(T)\|x\|_{{\cal D}(A)}.
\eeqn
\end{lemma}
\begin{proof}
First, from Lemma 2.3 in \cite{LT}, 
with the triplet $(L^p(\Om),A(t),f)$ being replaced by $(X,A,x)$, 
we have $\int_0^t{\rm e}^{\xi A}x\d\xi\in{\cal D}(A)$ for every $t>0$
and $A\int_0^t{\rm e}^{\xi A}x\d\xi={\rm e}^{tA}x-x$. Now, since $x\in{\cal D}(A)$
from the equality $A(\l-A)^{-1}x=(\l-A)^{-1}Ax$ and the definition 
of $\{{\rm e}^{tA}\}_{t>0}$ by Dunford integrals, we have 
$A{\rm e}^{tA}x={\rm e}^{tA}Ax$ for every $t>0$. Hence
\beqn\no
&&\hskip -1truecm
\int_0^t\|A{\rm e}^{\xi A}x\|_X\d\xi=\int_0^t\|{\rm e}^{\xi A}Ax\|_X\d\xi
\le\wtil c_0\int_0^t\xi^{(\b-1)/\a}\|Ax\|_X\d\xi
\le \wtil c_0\|x\|_{{\cal D}(A)}t^{(\a+\b-1)/\a}.
\eeqn
It follows that the map $\xi\to\|A{\rm e}^{\xi A}x\| _X$ belongs to $L^1((0,t);X)$
for all $t\in (0,T]$ and ${\rm e}^{tA}x-x=
A\int_0^t{\rm e}^{\xi A}x\d\xi=\int_0^tA{\rm e}^{\xi A}x\d\xi=\int_0^t{\rm e}^{\xi A}Ax\d\xi$.
Thus, (\ref{3.2}) with $n=0$ yields to
\beqn\label{4.10}
&&\hskip -1,3truecm
\|{\rm e}^{tA}x-x\|_{(X,{\cal D}(A))_{\g,p}}\le 
c_1(T)\int_0^t \xi^{(\b-1-\g)/\a}\|Ax\|_{X}\d\xi
\le c_1(T)\|x\|_{{\cal D}(A)}t^{(\a+\b-1-\g)/\a}.
\eeqn
Consequently, for every $t\in [0,T]$, we get
\beqn\label{4.11}
&&\hskip -1truecm
\|{\rm e}^{tA}x\|_{(X,{\cal D}(A))_{\g,p}}\le 
\|{\rm e}^{tA}x-x\|_{(X,{\cal D}(A))_{\g,p}}+\|x\|_{(X,{\cal D}(A))_{\g,p}}
\no
\\[2mm]
&&\hskip -1,4truecm
\le [c_1(T)t^{(\a+\b-1-\g)/\a}+c(\g,p)]\|x\|_{{\cal D}(A)}.
\eeqn
Now, let $\s\in(0,(\a+\b-1-\g)/\a)$ and $0<s<t\le T$
\footnote{Since $\s\in(0,(\a+\b-1-\g)/\a)$ and ${\rm e}^{0A}$ is defined to be $1$, 
the case $s=0$ follows from (\ref{4.10}).}. Then, reasoning as in the derivation of (\ref{4.10}), we obtain
\beqn\label{4.12}
\hskip -1,5truecm
\|{\rm e}^{tA}x-{\rm e}^{sA}x\|_{(X,{\cal D}(A))_{\g,p}}
\!\!\!&=&\!\!\!
\big\|\int_{s}^{t}A{\rm e}^{\xi A}x\d\xi\big\|_{(X,{\cal D}(A))_{\g,p}}\no
\\[1mm]
\hskip -1,5truecm
&\le&\!\!\! 
c_1(T)\|x\|_{{\cal D}(A)}(t-s)^{(\a+\b-1-\g)/\a}\no
\\[2mm]
\hskip -1,5truecm
&\le&\!\!\! 
c_1(T)T^{(\a+\b-1-\g-\a\s)/\a}\|x\|_{{\cal D}(A)}(t-s)^\s.
\eeqn
From (\ref{4.11}) and (\ref{4.12}) we deduce (\ref{4.9}) with 
\beqn\label{4.13}
&&\hskip -1truecm
c_4(T)=c(\g,p)+c_1(T)T^{(\a+\b-1-\g-\a\s)/\a}(T^\s+1)
\eeqn
 and the proof is complete.
\end{proof}
In the next section, in order to obtain optimal regularity
for solutions of degenerate parabolic equations, we will need to estimate
the maps $t\to {\rm e}^{tA}[f(t)-f(0)]$ and $t\to \int_0^tA{\rm e}^{(t-\xi)A}[f(\xi)-f(t)]\d\xi$,
where $f$ is H\"older continuous.
The next two lemmas give us the desired results in this direction.
\begin{lemma}\label{lem4.4}
Let $\a$, $\b\in(0,1]$ such that $\b<\a$ and $2\a+\b>2$. Then, for every
$\mu\in((2-\a-\b)/\a,1)$, $\g\in (0,\a\mu+\a+\b-2)$ and $\s\in (0,(\a\mu+\a+\b-2-\g)/\a)$
the linear operator 
\beqn\label{4.14}
&&\hskip -2truecm
[Q_2f](t):={\rm e}^{tA}[f(t)-f(0)]
\eeqn 
maps $C^\mu([0,T];X)$ into $C^\s([0,T];(X,{\cal D}(A))_{\g,p})$, $p\in[1,\infty]$, 
and satisfies the estimate:
\beqn\label{4.15}
&&\hskip -1truecm
\|Q_2f\|_{\s,T;(X,{\cal D}(A))_{\g,p}}
\le T^{(\a\mu+\a+\b-2-\g-\a\s)/\a}c_5(T)|f|_{\mu,T;X}.
\eeqn
\end{lemma}
\begin{proof}
When $f\in C^\mu([0,T];X)$, $\mu\in((2-\a-\b)/\a,1)$, 
from (\ref{3.2}) with $n=0$ it follows
\beqn\label{4.16}
&&\hskip -1truecm
\|[Q_2f](t)\|_{(X,{\cal D}(A))_{\g,p}}\le c_1(T)|f|_{\mu,T;X}t^{(\a\mu+\b-1-\g)/\a},
\q\forall\,t\in[0,T],
\eeqn
where $(\a\mu+\b-1-\g)/\a>0$, due to $0<\g<\a\mu+\a+\b-2\le\a\mu+\b-1$.
Now, let $\s\in (0,(\a\mu+\a+\b-2-\g)/\a)$ and $0<s<t\le T$
 \footnote{Since $0<\s<(\a\mu+\a+\b-2-\g)/\a\le (\a\mu+\b-1-\g)/\a$ the case $s=0$ follows from (\ref{4.16}).}. We have
\beqn\label{4.17}
&&\hskip -1,5truecm
\|[Q_2f](t)-[Q_2f](s)\|_{(X,{\cal D}(A))_{\g,p}}\le \sum_{k=1}^2I_k(s,t),
\eeqn
where 
\beqn\no
&&\hskip -1truecm
I_1(s,t):=\|{\rm e}^{tA}[f(t)-f(s)]\|_{(X,{\cal D}(A))_{\g,p}},\no
\\[2mm]
&&\hskip -1truecm
I_2(s,t):=\|({\rm e}^{tA}-{\rm e}^{sA})[f(s)-f(0)]\|_{(X,{\cal D}(A))_{\g,p}}.\no
\eeqn
Concerning $I_1(s,t)$, the same reasoning made to derive (\ref{4.16}) 
lead us to
\beqn\label{4.18}
&&\hskip -1,5truecm
I_1(s,t)\le c_1(T)t^{(\a\mu+\b-1-\g-\a\s)/\a}|f|_{\mu,T;X}(t-s)^\s,
\eeqn
the exponent $(\a\mu+\b-1-\g-\a\s)/\a$ being positive, since
$0<\s<(\a\mu+\a+\b-2-\g)/\a\le (\a\mu+\b-1-\g)/\a$. Instead, 
using (\ref{3.3}) with $n=0$ we obtain
\beqn\label{4.19}
&&\hskip -1,5truecm
I_2(s,t)\le \s^{-1}s^{(\a\mu+\a+\b-2-\g-\a\s)/\a}c_1(T)|f|_{\mu,T;X}(t-s)^\s
\eeqn
Therefore, (\ref{4.16})--(\ref{4.19}) yield to (\ref{4.15}) with 
\beqn\label{4.20}
&&\hskip -1,5truecm
c_5(T)=c_1(T)[T^{(1-\a)/\a}(T^\s+1)+\s^{-1}].
\eeqn
The proof is now complete.
 \end{proof}
\begin{lemma}\label{lem4.5}
Let $\a$, $\b\in(0,1]$ such that $\b<\a$ and $3\a+\b>3$.
Then, for every  $\mu\in((3-2\a-\b)/\a,1)$, $\g\in (0,\a\mu+2\a+\b-3)$ and
$\s\in (0,(\a\mu+2\a+\b-3-\g)/\a)$ the linear operator
\beqn\label{4.21}
&&\hskip -2truecm
[Q_3f](t):=\int_0^tA{\rm e}^{(t-\xi)A}[f(\xi)-f(t)]\d\xi 
\eeqn
maps $C^\mu([0,T];X)$ into $C^\s([0,T];(X,{\cal D}(A))_{\g,p})$, $p\in[1,\infty]$, 
and satisfies the estimate:
\beqn\label{4.22}
&&\hskip -1truecm
\|Q_3f\|_{\s,T;(X,{\cal D}(A))_{\g,p}}
\le T^{(\a\mu+2\a+\b-3-\g-\a\s)/\a}c_6(T)|f|_{\mu,T;X}.
\eeqn
\end{lemma}
\begin{proof}
Let $f\in C^\mu([0,T];X)$, $\mu\in((3-2\a-\b)/\a,1)$. 
Then, for every $t\in[0,T]$, using (\ref{3.2}) with $n=1$ we find
\beqn\label{4.23}
\hskip -1,5truecm
\|[Q_3f](t)\|_{(X,{\cal D}(A))_{\g,p}}
\!\!\!&\le&\!\!\!
c_1(T)|f|_{\mu,T;X}\int_0^t(t-\xi)^{(\a\mu+\b-2-\g)/\a}\d\xi\no
\\[1mm]
\hskip 1,5truecm
&\le&\!\!\!\frac{\a c_1(T)}{\a\mu+\a+\b-2-\g}|f|_{\mu,T;X}t^{(\a\mu+\a+\b-2-\g)/\a}.
\eeqn
Notice that the choice $\g\in (0,\a\mu+2\a+\b-3)$ implies
$(\a\mu+\a+\b-2-\g)/\a>0$ in the latter inequality.
Now, let $\s\in (0,(\a\mu+2\a+\b-3-\g)/\a)$ and $0<s<t\le T$
\footnote{Since $\s\in (0,(\a\mu+2\a+\b-3-\g)/\a)$, 
the case $s=0$ follows from (\ref{4.23}) once we observe that
$$(\a\mu+\a+\b-2-\g)/\a-\s
>\big(\a\mu+\a+\b-2-\g-(\a\mu+2\a+\b-3-\g)\big)/\a=(1-\a)/\a\ge 0.$$.}. 
We have
\beqn\label{4.24}
&&\hskip -1,5truecm
\|[Q_3f](t)-[Q_3f](s)\|_{(X,{\cal D}(A))_{\g,p}}\le \sum_{k=1}^3J_k(s,t),
\eeqn
where
\beqn
&&\hskip -1truecm
J_1(s,t):=\big\|\int_0^{s}[A{\rm e}^{(t-\xi)A}-A{\rm e}^{(s-\xi)A}]
[f(\xi)-f(s)]\d\xi\big\|_{(X,{\cal D}(A))_{\g,p}},\no
\\[1mm]
&&\hskip -1truecm
J_2(s,t):=\big\|\int_0^{s}A{\rm e}^{(t-\xi)A}
[f(s)-f(t)]\d\xi\big\|_{(X,{\cal D}(A))_{\g,p}},\no
\\[1mm]
&&\hskip -1truecm
J_3(s,t):=\big\|\int_{s}^{t}A{\rm e}^{(t-\xi)A}
[f(\xi)-f(t)]\d\xi\big\|_{(X,{\cal D}(A))_{\g,p}}.\no
\eeqn
We examine first $J_1(s,t)$. To this purpose, from inequality (\ref{3.3})
with $n=1$ we deduce
\beqn\label{4.25}
\hskip -1truecm
J_1(s,t)\!\!\!&\le&\!\!\!
\s^{-1}c_1(T)|f|_{\mu,T;X}
\big[\int_0^{s}(s-\xi)^{(\a\mu+\a+\b-3-\g-\a\s)/\a}\d\xi\big](t-s)^\s\no
\\[2mm]
\hskip -1truecm
&\le&\!\!\!
\frac{\s^{-1}\a c_1(T)}{\a\mu+2\a+\b-3-\g-\a\s}|f|_{\mu,T;X}
s^{(\a\mu+2\a+\b-3-\g-\a\s)/\a}(t-s)^\s.
\eeqn
Let us turn to $J_2(s,t)$. Since $\a+\b-2-\g<0$,
inequality (\ref{3.2}) with $n=1$ yields to
\beqn\label{4.26}
\hskip -0,5truecm
J_2(s,t)\!\!\!&\le&\!\!\! 
c_1(T)|f|_{\mu,T;X}\big[\int_0^{s}(t-\xi)^{(\b-2-\g)/\a}\d\xi\big](t-s)^{\mu}\no
\\[1mm]
\hskip -0,5truecm
&\le&\!\!\!
\frac{\a c_1(T)}{2+\g-\a-\b}|f|_{\mu,T;X}
[(t-s)^{(\a+\b-2-\g)/\a}-t^{(\a+\b-2-\g)/\a}](t-s)^{\mu}\no
\\[2mm]
\hskip -0,5truecm
&\le&\!\!\!
\frac{\a c_1(T)}{2+\g-\a-\b}|f|_{\mu,T;X}(t-s)^{(\a\mu+\a+\b-2-\g)/\a}\no
\\[2mm]
\hskip -0,5truecm
&\le&\!\!\!
\frac{\a c_1(T)}{2+\g-\a-\b}|f|_{\mu,T;X}
t^{(\a\mu+\a+\b-2-\g-\a\s)/\a}(t-s)^\s.
\eeqn
Finally, concerning $J_3(s,t)$, still from (\ref{3.2}) with $n=1$ we get
\beqn\label{4.27}
\hskip -1,5truecm
J_3(s,t)\!\!\!&\le&\!\!\! 
c_1(T)|f|_{\mu,T;X}\int_{s}^{t}(t-\xi)^{(\a\mu+\b-2-\g)/\a}\d\xi\no
\\[1mm]
\hskip -1,5truecm
&\le&\!\!\!
\frac{\a c_1(T)}{\a\mu+\a+\b-2-\g}|f|_{\mu,T;X}(t-s)^{(\a\mu+\a+\b-2-\g)/\a}\no
\\[2mm]
\hskip -1,5truecm
&\le&\!\!\!
\frac{\a c_1(T)}{\a\mu+\a+\b-2-\g}|f|_{\mu,T;X}
t^{(\a\mu+\a+\b-2-\g-\a\s)/\a}(t-s)^\s.
\eeqn
As a consequence, replacing (\ref{4.25})--(\ref{4.27}) in (\ref{4.24}), we obtain
\beqn\label{4.28}
&&\hskip -1,5truecm
\|[Q_3f](t)-[Q_3f](s)\|_{(X,{\cal D}(A))_{\g,p}}\le
T^{(\a\mu+2\a+\b-3-\g-\a\s)/\a}c_7(T)|f|_{\mu,T;X}(t-s)^\s,
\eeqn
where $c_7(T)=\a c_1(T)c_8(T)$, $c_8(T)$ being defined by
\beqn
&&\hskip -1truecm
c_8(T)=
\Big[\frac{\s^{-1}}{\a\mu+2\a+\b-3-\g-\a\s}
+\frac{\a\mu T^{(1-\a)/\a}}{(2+\g-\a-\b)(\a\mu+\a+\b-2-\g)}\Big].\no
\eeqn
Summing up (\ref{4.23}) and (\ref{4.28}) we easily derive (\ref{4.22}) with 
\beqn\label{4.29}
&&\hskip -2truecm
c_6(T)=\a c_1(T)\big[(\a\mu+\a+\b-2-\g)^{-1}T^{(1-\a+\a\s)/\a}+c_8(T)\big].
\eeqn
This completes the proof.
\end{proof}
\section{Maximal regularity for degenerate equations}\label{Sec5}
\setcounter{equation}{0}
In this section we apply the preliminary lemmata of Section \ref{Sec4} 
for proving time and space regularity of solutions to the following 
degenerate first-order initial value problem
\beqn\label{5.1}
&&\hskip -2truecm
\left\{\!\!
\begin{array}{lll}
D_t(Mv(t))=Lv(t)+f(t), \qq t\in(0,T],
\\[2mm]
Mv(0)=u_0,
\end{array}
\right.
\eeqn
in a Banach space $X$. Here  $M$ and $L$ are two closed linear operators 
from $X$ to itself having domains, respectively, 
${\cal D}(M)$ and ${\cal D}(L)$ with ${\cal D}(L)\subset {\cal D}(M)$, 
$f\in C([0,T];X)$ is a given function and $u_0\in X$ is a given initial value.

We want to stress that, since the Cauchy problem (\ref{5.1}) 
coincides with that of type (D-E.1) in \cite[Section 3.3]{FY1},
the range of the possible applications of our results turns out to be very large.
To this purpose, we refer the interested reader to \cite{FY1}, and to
the references therein, for a list of boundary value problems related to degenerate parabolic equations which can be reduced to (\ref{5.1}) via an abstract reformulation.

According to \cite{FY1}, we recall that the 
$M$-modified resolvent set of $L$ is the set
$\rho_M(L)=\{\l\in\csp:\l M-L$ {\rm has ``bounded inverse"}
$M(\l M-L)^{-1}$ {\rm on} $X\}$. It is easy to prove that 
$\rho_M(L)\subset\rho(LM^{-1})$ and that
$M(\l M-L)^{-1}=(\l-LM^{-1})^{-1}$, $\l\in\rho_M(L)$
(cf. \cite[Theorem 1.14]{FY1}).
With the notion of $M$-modified resolvent set of $L$ at hand, we assume:
\begin{itemize}
\item[(H2)] $\rho_M(L)$ contains the region
$\S=\{\l\in\csp:\Re{\rm e}\l\ge -c(|\Im{\rm m}\l|+1)^\a\}$ and,
for every $\l\in\S$, the following estimate holds
\beqn
&&\hskip -2truecm
\|M(\l M-L)^{-1}\|_{{\cal L}(X)}\le C(|\l|+1)^{-\b},\no
\eeqn 
for some exponents $0<\b<\a\le 1$ and constants $c$, $C>0$.
\end{itemize}
Of course, assumption (H2) implies that the operator $A=LM^{-1}$ with domain 
${\cal D}(A)=M({\cal D}(L))$ satisfies assumption (H1) and hence that
it generates a semigroup $\{{\rm e}^{tA}\}_{t\ge 0}$ 
defined by (\ref{2.1}) and satisfying (\ref{2.2}). 

Notice that, due to the identity $L(\l M-L)^{-1}=\l M(\l M-L)^{-1}-I$, (H2) 
reads equivalently to
\beqn\label{5.2}
&&\hskip -2truecm
\|L(\l M-L)^{-1}\|_{{\cal L}(X)}\le C|\l|(|\l|+1)^{-\b}+1\le (C+1)(|\l|+1)^{1-\b}.
\eeqn 
However, until now, under assumption (\ref{5.2}) 
{\it only} results of time regularity have been established.
See, for instance, \cite[Theorem 9]{FY2} and \cite[Theorem 7.2]{FLT}.
A result of space regularity has been obtained in \cite{FY2},
but with a stronger hypothesis of abstract potential
type on the operator $T=ML^{-1}=A^{-1}$, precisely
\beqn\label{5.3}
\|L(\l M-L)^{-1}\|_{{\cal L}(X)}=\|(\l T-I)^{-1}\|_{{\cal L}(X)}\le C, \qq \l>0.
\eeqn  
In this case ${\cal T}$, 
the part of $T$ in the closure $\ov{R(T)}$ of its range, 
has a densely defined inverse ${\cal T}^{-1}$ (unbounded, in general) which
generates an analytic semigroup in $\ov{R(T)}$. 
Then, denoted by $P$ the projection operator onto 
the null space $N(T)$ of $T$ and provided
that some suitable assumptions are satisfied on 
$(I-P)f$ and $(I-P)Lu_0$, $v_0=Mu_0$,
in \cite[Theorem 5]{FY2} it is shown that 
$D_tMv$ belongs to $B([0,T],X_{{\cal T}^{-1}}^{\theta,\infty})$, $\theta\in(0,1)$. 
This is done by means of customary techniques of analytic semigroup theory.
In particular, since assumption (\ref{5.3}) implies
the identity $X_{{\cal T}^{-1}}^{\theta,\infty}=
(\ov{R(T)},{\cal D}({\cal T}^{-1}))_{\theta,\infty}$, 
the quoted result extends \cite[Theorem 5.5]{SI} to degenerate equations.

The main problem in \cite{FY2} lies in the characterization of projection $P$,
which is crucial  when space regularity is investigated. From this point of view, our aim
is twofold. At first, to replace (\ref{5.3}) with the more general assumption (H2)
removing the analyticity of the semigroup ${\rm e}^{tA}$. Then, to show both
time and space regularity for $D_tMv$ without invoking $P$.

We begin proving two theorems concerning the regularity of the strict solution to (\ref{5.1}). 
Recall that, according to \cite[page 53]{FY1}, by a strict solution $v$ to (\ref{5.1})
we mean a function $v\in C((0,T],{\cal D}(L))$ such that
$Mv\in C^1((0,T];X)$ and (\ref{5.1}) holds, where $Mv(0)=u_0$ is understood 
in the sense that $\lim_{t\to 0}\|ML^{-1}(Mv(t)-u_0)\|_X=0$.
\begin{theorem}\label{thm5.1}
Let assumption {\rm (H2)} be fulfilled with $2\a+\b>2$ and let $u_0\in M({\cal D}(L))$
and $f\in C^\mu([0,T];X)$, $\mu\in((2-\a-\b)/\a,1)$. 
Then, for every $\g\in (0, 2\a+\b-2)$, $\s\in (0,(2\a+\b-2-\g)/\a)$ and $p\in[1,\infty]$, problem $(\ref{5.1})$ has a unique strict solution $v$ such that
\beqn\label{5.4}
&&\hskip -2truecm
Mv\in C^1((0,T];X)\cap C([0,T];X)\cap C^{\s}([0,T];(X,{\cal D}(A))_{\g,p}).
\eeqn 
Moreover, the following estimate holds true:
\beqn\label{5.5}
&&\hskip -1,5truecm
\|Mv\|_{\s,T;(X,{\cal D}(A))_{\g,p}}
\le c_4(T)\|u_0\|_{{\cal D}(A)}+T^{\nu}c_{2}(T)\|f\|_{0,T;X},
\eeqn
where $\nu=(2\a+\b-2-\g-\a\s)/\a$. Here
$c_2(T)$ and $c_4(T)$ are the positive nondecreasing functions of $T$ defined, respectively, in $(\ref{4.5})$ and $(\ref{4.13})$.
\end{theorem}
\begin{proof}
First, when $f\in C^{\mu}([0,T];X)$, $\mu\in ((2-\a-\b)/\a,1)$, 
and $u_0\in M({\cal D}(L))= {\cal D}(A)$,
\cite[Theorem 3.9]{FY1} ensures that problem (\ref{5.1}) 
admits a unique strict solution $v$ such that 
$Mv\in C^1((0,T];X)\cap C([0,T];X)$. 
In particular, the following representation holds:
\beqn\label{5.6}
&&\hskip -1,5truecm
(Mv)(t)={\rm e}^{tA}u_0+[Q_1f](t), \qq t\in[0,T],
\eeqn
$Q_1$ being defined in (\ref{4.1}). This is a consequence of  
\cite[Theorem 3.7]{FY1} and the Remark to it,
changing the unknown function to $w=Mv$ and rewriting 
(\ref{5.1}) into the equivalent form
\beqn\no
&&\hskip -1,5truecm
D_tw(t)=Aw(t)+f(t),\;\ t\in(0,T],\q w(0)=u_0.
\eeqn
Further, for every $\g\in (0, 2\a+\b-2)$ and $\s \in (0,(2\a+\b-2-\g)/\a)$,
Lemmas \ref{lem4.1} and \ref{lem4.3} imply that 
$Q_1f$ and ${\rm e}^{\cdot A}u_0$ belong to
$C^{\s}([0,T];(X,{\cal D}(A))_{\g,p})$, $p\in[1,\infty]$, and the same
assert is true for $Mv$ by virtue of (\ref{5.6}). 
Finally, estimate (\ref{5.5}) follows from (\ref{4.2}), (\ref{4.9}) and (\ref{5.6}).
\end{proof}
\begin{remark}\label{rem5.2}
\emph{We stress that, even if in Theorem \ref{thm5.1} $f$ is assumed 
$\mu$-H\" older continuous in time,  during the proof we have used 
Lemma \ref{lem4.1} which requires only the mere continuity of $f$. 
This is, for, $\mu$ being in $((2-\a-\b)/\a,1)$, it is not guaranteed 
that $f\in C^{\s}([0,T];X)$, $\s\in (0,(\a+\b-1-\g)/\a)$, $\g\in(0,\a+\b-1)$, 
in order to apply Lemma \ref{lem4.2}. 
Indeed, provided $\a+\b>3/2$ and $\g\in (0, 2(\a+\b)-3)\subset(0,\a+\b-1)$,
it may happen that $(2-\a-\b)/\a<\mu<\s<(\a+\b-1-\g)/\a$, 
so that $f\notin C^{\s}([0,T];X)$. 
Since $\a+\b>3/2$ implies $2\a+\b>2$, such a case may effectively take place
if $\a$ and $\b$ are large enough. 
Situation is not better if we try to apply Lemma \ref{lem4.2} 
restricting $\g$ and $\s$ to vary in the smaller intervals of Lemma \ref{lem4.1}. 
In fact, it may occur that $(2-\a-\b)/\a<\mu<\s<(2\a+\b-2-\g)/\a$, provided 
$3\a+2\b>4$ and $\g\in (0,3\a+2\b-4)\subset(0,2\a+\b-2)$.}
\end{remark}
\begin{theorem}\label{thm5.3}
Let  the assumptions in {\rm Theorem \ref{thm5.1}} be satisfied, but with $\a+\b>3/2$. 
Then, for every $\g\in [2(\a+\b)-3, \a+\b-1)$, $\s \in (0,(\a+\b-1-\g)/\a)$ 
and $p\in[1,\infty]$, problem $(\ref{5.1})$ has a unique strict solution $v$ 
satisfying $(\ref{5.4})$. 
Moreover, the following estimate holds true:
\beqn\label{5.7}
&&\hskip -1truecm
\|Mv\|_{\s,T;(X,{\cal D}(A))_{\g,p}}
\le c_4(T)\|u_0\|_{{\cal D}(A)}+T^{\nu}c_3(T)\max\{1,T^{\mu-\s}\}\|f\|_{\mu,T;X},
\eeqn
where $\nu=(\a+\b-1-\g-\a\s)/\a$. 
Here $c_3(T)$ is the positive nondecreasing function of $T$ defined in $(\ref{4.8})$.
\end{theorem}
\begin{proof}
As before, since $\a+\b>3/2$ implies $2\a+\b>2$, 
the belonging of $Mv$ to $C^1((0,T];X)\cap C([0,T];X)$ follows 
from \cite[Theorem 3.9]{FY1}. Now, our assumptions on $u_0$, $\g$ and $\s$
enable us to use Lemma \ref{lem4.3} to ensure that  
${\rm e}^{\cdot A}u_0$ belong to $C^{\s}([0,T];(X,{\cal D}(A))_{\g,p})$, $p\in[1,\infty]$. 
In addition, the assumption $\g\ge2(\a+\b)-3$ imply the following chain of inequalities
\beqn\no
&&\hskip -2truecm
0<\s<(\a+\b-1-\g)/\a\le (2-\a-\b)/\a<\mu<1.
\eeqn
As a consequence, we have $f\in C^\s([0,T];X)$, $\s\in (0,(\a+\b-1-\g)/\a)$, 
and we are in position to apply Lemma \ref{lem4.2}. 
Thus, $Q_1f$ is in $C^{\s}([0,T];(X,{\cal D}(A))_{\g,p})$, $p\in[1,\infty]$, 
and satisfies (\ref{4.6}), with $g$ being replaced by $f$.
Hence, from (\ref{5.6}) we deduce that 
$Mv\in C^{\s}([0,T];(X,{\cal D}(A))_{\g,p})$, $p\in[1,\infty]$.
Finally, estimate (\ref{5.7}) follows from (\ref{4.6}), (\ref{4.9}), (\ref{5.6})
and the inequality $\|f\|_{\s,T;X}\le\max\{1,T^{\mu-\s}\}\|f\|_{\mu,T;X}$.
\end{proof}
We now come to our main theorem, which provides regularity in both time and
space for the derivative $D_tMv$. The following statement improves
\cite[Theorem 3.26]{FY1} and the results in \cite{FY2} and \cite{FLT} mentioned before.
\begin{theorem}\label{thm5.4}
Let assumption {\rm (H2)} be fulfilled with $3\a+\b>3$ and let $u_0\in M({\cal D}(L))$
and $f\in C^\mu([0,T];X)$, $\mu\in((3-2\a-\b)/\a,1)$. Further, let assume that
\beqn\label{5.8}
&&\hskip -2truecm
Lv_0+f(0)=:g_0\in M({\cal D}(L)),\q u_0=Mv_0,\q v_0\in{\cal D}(L).
\eeqn
Then, for every $\g\in (0,\a\mu+2\a+\b-3)$, $\s\in (0,(\a\mu+2\a+\b-3-\g)/\a)$ 
and $p\in[1,\infty]$, problem $(\ref{5.1})$ has a unique strict solution $v$ such that
\beqn\no
&&\hskip -2truecm
D_tMv\in C^\s([0,T];(X,{\cal D}(A))_{\g,p}).
\eeqn
Moreover, the following estimate holds true:
\beqn\label{5.9}
&&\hskip -2truecm
\|D_tMv\|_{\s,T;(X,{\cal D}(A))_{\g,p}}
\le c_4(T)\|g_0\|_{{\cal D}(A)}+T^\nu C(T)|f|_{\mu,T;X},
\eeqn
where $\nu=(\a\mu+2\a+\b-3-\g-\a\s)/\a$ and $C(T)=T^{(1-\a)/\a}c_5(T)+c_6(T)$. 
Here $c_5(T)$ and $c_6(T)$ are the positive nondecreasing functions of $T$ 
defined, respectively, in $(\ref{4.20})$ and $(\ref{4.29})$.
\end{theorem}
\begin{proof}
First, since $3\a+\b>3$ implies $2\a+\b>2$ and 
$f\in C^{\mu}([0,T];X)$, where $\mu\in((3-2\a-\b)/\a,1)\subset((2-\a-\b)/\a,1)$,
Theorem \ref{thm5.1} applies and (\ref{5.6}) holds.
Hence, as shown in \cite[Remark page 55]{FY1}, 
differentiating (\ref{5.6}) with respect to $t$ and using (\ref{5.8}), we deduce
\beqn\label{5.10}
&&\hskip -2truecm
D_t(Mv(t))={\rm e}^{tA}g_0+[Q_2f](t)+[Q_3f](t),\q t\in[0,T],
\eeqn
the $Q_j$'s, $j=2,3$, being defined, respectively, in (\ref{4.14}) and (\ref{4.21}).
In particular, $Mv\in C^1([0,T];X)$ and
the equation in (\ref{5.1}) makes sense even at $t=0$. 
Now, notice that
\beqn\no
&&\hskip -1truecm
\g\in (0,\a\mu+2\a+\b-3)\subset(0,\a\mu+\a+\b-2)\subset(0,\a+\b-1),
\\[2mm]
&&\hskip -1truecm
\s\in \Big(0,\frac{\a\mu+2\a+\b-3-\g}{\a}\Big)\subset 
\Big(0,\frac{\a\mu+\a+\b-2-\g}{\a}\Big)\subset\Big(0,\frac{\a+\b-1-\g}{\a}\Big),\no
\eeqn
so that all the assumptions of Lemmas \ref{lem4.3}--\ref{lem4.5} are satisfied. 
Therefore, for every $\g\in (0,\a\mu+2\a+\b-3)$
and $\s\in (0,(\a\mu+2\a+\b-3-\g)/\a)$,
${\rm e}^{\cdot A}g_0$ and $Q_jf$, $j=2,3$,
belong to $C^\s([0,T];(X,{\cal D}(A))_{\g,p})$, $p\in[1,\infty]$.
Of course, due to (\ref{5.10}), the same belonging holds for the derivative $D_tMv$. 
Finally, estimates (\ref{4.9}), (\ref{4.14}) and (\ref{4.22}) yields to (\ref{5.9}),
and the proof is complete.
\end{proof}
\centerline\acks
The author is really grateful to Professor Angelo Favini of the Universit\`a di Bologna
for having him suggested the reference in \cite{WI}.

\end{document}